\documentclass[runningheads]{llncs}
\usepackage{graphicx}

\usepackage{amssymb,amsmath}
\usepackage{xcolor}
\usepackage{hyperref}
\usepackage[all]{xy}

\begin{document}
\title{On the functorial properties of the $p$-analog of the Fourier-Stieltjes algebras}
\titlerunning{functorial properties of $B_p(G)$}
%

\author{Mohammad Ali Ahmadpoor\inst{1}\orcidID{0000-0001-6902-1916} \and Marzieh Shams Yousefi \inst{2}\orcidID{0000-0003-0426-708X}}
\authorrunning{M. A. Ahmadpoor \& M. Shams Yousefi}
%
\institute{${}^{{1,2}}$ Department of Pure Mathematics, Faculty of Mathematical Sciences, University of Guilan, Rasht, Iran\\
\email{m-a-ahmadpoor@phd.guilan.ac.ir}\\
\email{m.shams@guilan.ac.ir}}

\maketitle              
\begin{abstract}
In this paper, some known results will be generalized. Firstly, the idempotent theorem on the Fourier-Stieltjes algebra will be promoted and linked to the $p$-analog of such an algebra. Next, the  $p$-analog of the $\pi$-Fourier space introduced by Arsac will be given, and by taking advantage of the theory of ultra filters, the connection between the dual space of the algebra of $p$-pseudofunctions and the $p$-analog of the $\pi$-Fourier space, will be fully investigated.

As the main result, one of the significant and applicable functorial properties of the $p$-analog of the Fourier-Stieltjes algebras will be achieved. 

\end{abstract}

\section{Introduction}\label{subsection1.1}

For a locally compact group $G$, the Fourier algebra, $A(G)$, and the Fourier-Stieltjes algebra, $B(G)$, have been found by Eymard in 1964 \cite{EYMARD1964}. He investigated almost all functorial properties of such algebras. On the other hand, idempotent elements of $B(G)$ is introduced by Host \cite{HOST1986}, and has gotten accurate by Ilie and Spronk in \cite{ILIESPRONK2005}. Even Runde \cite{RUNDE2007} went beyond and add some specific conditions to it, by benefiting from the theory of uniformly convex Banach spaces. In the next attempt on studying Fourier type algebras, for a representation $(\pi,\mathcal{H})$ of $G$, on a Hilbert space $\mathcal{H}$, Arsac \cite{ARSAC1976} introduced $\pi$-Fourier and $\pi$-Fourier-Stieltjes spaces $A_\pi$ and $B_\pi$, and tremendously studied their functorial properties. Meanwhile, Fig\`a-Talamanca-Herz algebra was initially defined for abelian locally compact groups \cite{FIGATALAMANCA1965}, and then for general locally compact groups in \cite{HERZ1971}. Afterwards, in \cite{RUNDE2005}, Runde took main step and determined the true $p$-analog of Eymard's $B(G)$, as it is indicated via $B_p(G)$. He has indicated that the space $B_p(G)$ is a communicative unital Banach algebra and in the case that the underlying group $G$ is amenable, it can be identified with the multiplier algebra of the Fig\`a-Talamanca-Herz algebra, i.e. the Banach space $\mathcal{M}(A_p(G))$. Accordingly, the vast majority of functorial properties of the $p$-analog of the Fourier-Stieltjes algebras $B_p(G)$ have been remained unknown, and they would be the cause for a huge amount of studies. For instance, in \cite{NEUFANGRUNDE2007}, one of the possible $p$-operator space structure on $B_p(G)$ is studied while another one is introduced by authors \cite{AHSH12020}. Besides, there can be found considerable amount of primary questions about the element of such algebras is still open. In this paper, we try to come up with detailed explanation, to touch on some topics, which can be considered as initial steps in the studying on the $p$-analog of the Fourier-Stieltjes algebras. Current paper is organized as following: In Section \ref{SECTIONPRELIMINARIES}, we give some preliminaries about representations of a locally compact group on $QSL_p$s-spaces which are building blocks of Runde's $B_p(G)$. It is indicated that the communicative Banach space $B_p(G)$ is the dual space of the algebra of universal $p$-pseudofunctions $UPF_{p}(G)$. Next, in the main section, Section \ref{SECTIONGENERALIZATIONS}, we divide our results into three main subsections. In one approach, the Subsection \ref{subsectionIT} is devoted to the generalization of the idempotent theorem (Theorem \ref{IDMTHEOREM}) in terms of $p$-analog of the Fourier-Stieltjes algebra. Indeed, we have benefited from the aforementioned theory of Banach spaces and added one more equivalent statement to the latest version of it, that is done by Runde \cite[Theorem 1.5]{RUNDE2007}. In the next attempt, we have touched briefly on the generalization of $\pi$-Fourier spaces, and simply introduced it, then, through restating the result in \cite{RUNDE2005}, in respect to the notion of ultrafilter, we have reached to the appropriate description of the dual of the algebra of $p$-pseudofunctions $PF_{p,\pi}(G)$, for an arbitrary representation $(\pi,E)$ of the locally compact group $G$ on a $QSL_p$-space $E$, that is denoted by $B_{p,\pi}$ (Proposition \ref{PROPULTRAAP}). At this aim, we have generated somewhat crucial properties on the representations of a locally compact group, via the notion of ultrafilter and ultrapower space.\\
Finally, as a conclusion of previous results, in Proposition \ref{PROPEXTENSION} we have shown that a function in $u\in B_p(G_0)$ can be extended to a function $u^\circ\in B_p(G)$, where $G_0$ is an open amenable subgroup of the locally compact group $G$, and the general form of this proposition channels us to Theorem \ref{THEOREMCONCLUSION}.

\section{Preliminaries}\label{SECTIONPRELIMINARIES}
In this paper, $G$ and $H$ are locally compact groups, and for $p\in (1,\infty)$, the number $p'$ is its complex conjugate, i.e. $1/p+1/p'=1$. In the first step, we give essential notions and definitions on $QSL_p$-spaces, and representations of groups on such spaces. For more information one can see \cite{RUNDE2005}.
\begin{definition}\label{def1}
A representation of a locally compact group $ G $ is a
pair $ (\pi , E) $, where $ E $ is a Banach space and $ \pi $ is a group homomorphism from $ G $ into the
invertible isometries on $ E $, that is continuous with respect to the given topology on $ G $
and the strong operator topology on $ \mathcal{B}(E) $.
\end{definition}

\begin{remark}
Every representation $ (\pi , E) $ of a locally compact group $ G $ induces a representation of the group algebra $L_1(G)$ on $ E $, i.e. a contractive algebra homomorphism from $L_1(G)$ into
$ \mathcal{B}(E) $, which we shall denote likewise by $ \pi $, through
\begin{align}\label{11}
&\pi(f)=\int f(x)\pi (x)dx,\ f\in L_1(G),\\
&\langle \pi(f)\xi, \eta\rangle =\int f(x)\langle \pi (x)\xi, \eta\rangle dx,\quad \xi\in E,\:\eta\in E^*,\nonumber
\end{align}
where the integral \eqref{11} converges with respect to the strong operator topology.
\end{remark}

\begin{definition}
Let $ (\pi , E)$ and $(\rho , F) $ be representations of the locally compact group $G$. Then

\begin{enumerate}
\item
$ (\pi , E)$ and $ (\rho , F) $ are called equivalent, if there exists an invertible isometry $\varphi : E\rightarrow F$ such that
\begin{equation*}
\varphi\pi (x)\varphi^{-1}=\rho (x),\qquad x\in G.
\end{equation*}
\item
$ (\rho , F) $ is said to be a subrepresentation of $ (\pi , E)$, if $F$ is a closed subspace of $E$, and for every $x\in G$ we have $\pi(x){|_{F}}=\rho(x)$.
\item
$ (\rho , F) $ is said to be contained in $ (\pi , E)$, if it is equivalent to a subrepresentation of $  (\pi , E) $, and will be denoted by $(\rho , F) \subset (\pi , E)$.
\end{enumerate}
\end{definition}

\begin{definition}
\begin{enumerate}
\item A Banach space is called an $L_p$-space if it is of the form $L_p(X)$ for some measure
space $X$.
\item A Banach space is called a $ QSL_p $-space if it is isometrically isomorphic to a quotient of a subspace of an $ L_p $-space.
\end{enumerate}
\end{definition}
We denote by $ \text{Rep}_p(G) $ the collection of all (equivalence classes) of representations of $G $ on a  $ QSL_p $-space.

\begin{definition}
A representation of a Banach algebra $\mathcal{A}$ is a pair $ (\pi , E)$, where $E$ is a Banach space, and $\pi $ is a contractive algebra homomorphism from $\mathcal{A}$ to $\mathcal{B}(E)$. We call
$ (\pi , E)$ isometric if $\pi $ is an isometry and essential if the linear span of $\{\pi(a)\xi \: : \ a \in \mathcal{A},\: \xi\in
E\} $ is dense in $E$. 
\end{definition}

\begin{remark}\label{rem11}
If $ G $ is a locally compact group and $ (\pi , E)$ is a representation of $ G $ in the sense of Definition
\ref{def1}, then \eqref{11} induces an essential representation of $L_1(G)$. Conversely,
every essential representation of $L_1(G)$ arises in this fashion.
\end{remark}

\begin{definition}

\begin{enumerate}
\item
A representation $ (\pi , E)\in \text{Rep}_p(G) $ is called cyclic, if there exists $\xi_0\in E$ such that $\pi (L_1(G))\xi_0$ is dense in $E$. The set of cyclic representations of group $G$ on $QSL_p$-spaces is denoted by $\text{Cyc}_p(G)$.
\item
A representation $ (\pi , E)\in \text{Rep}_p(G)$ is called $p$-universal, if it contains every cyclic representation.
\end{enumerate}

\end{definition}

\begin{remark}\label{REMARKGARDELLA}
By  \cite[Remark 2.9-$(3)$]{GARDELLATHIEL2015}, and \cite[Proposition 2.4]{GARDELLATHIEL2015}, it is easy to see that every $p$-universal representation of $G$, contains every cyclic representation of $G$ on a $QSL_p$-space, in the sense of equivalency. In Addition, every representation in $\text{Rep}_p(G)$ is contained in a $p$-universal representation. Actually, one could make a new $p$-universal representation by constructing $l_p$-direct sum of an arbitrary representation with a $p$-universal representation.
\end{remark}

Now we are ready to describe the Fig\`a-Talamanca-Herz, and the $p$-analog of the Fourier-Stieltjes algebras.

\begin{definition}
Fig\`a-Talamanca-Herz algebra on the locally compact group $G$, which is denoted by $A_p(G)$, is the collection of functions $u:G\rightarrow\mathbb{C}$ of the form

\begin{equation}\label{Ap}
u(x)=\sum_{n=1}^{\infty}\langle \lambda_{p}(x)\xi_n,\eta_n\rangle,\quad x\in G,
\end{equation}
with 

\begin{equation}\label{Ap2}
(\xi_n)_{n\in\mathbb{N}}\subset L_{p}(G),\quad (\eta_n)_{n\in\mathbb{N}}\subset L_{p'}(G),\quad\text{and}\quad\sum_{n=1}^\infty\|\xi_n\| \|\eta_n\|<\infty,
\end{equation}
where $ \lambda_{p}$ is the left regular representation of $G$ on $L^{p}(G)$, defined as
\begin{align*}
&\lambda_{p}: G\rightarrow \mathcal{B}(L^{p}(G)),\quad\lambda_{p}(x)\xi(y)=\xi(x^{-1}y),\quad \xi\in L^{p}(G),\:x,y\in G.
\end{align*}
The norm of $A_p(G)$ is defined as 
\begin{equation*}
\| u\|=\inf\Big\{\sum_{n=1}^\infty\|\xi_n\| \|\eta_n\|\: :\ u(\cdot)=\sum_{n=1}^\infty\langle \lambda_{p}(\cdot)\xi_n, \eta_n\rangle\Big\},
\end{equation*}
where the infimum is taken over all expressions of $ u $ in \eqref{Ap} with \eqref{Ap2}.
With this norm and pointwise operations, $A_p(G)$ turns into a commutative regular Banach algebra.
\end{definition}

\begin{remark}
 The $p$-analog of the Fourier-Stieltjes algebra has been studied, for example in \cite{COWLING1979}, \cite{FORREST1994}, \cite{MIAO1996} and \cite{PIER1984}, as the multiplier algebra of the Fig\`a-Talamanca-Herz algebra. In this paper, we follow the construction of Runde in definition and notation (See \cite{RUNDE2005}) which we swap indexes $p$ and $p'$.
\end{remark}

\begin{definition}\label{definitionofBp}
The set of all functions of the form
\begin{equation}\label{eqq1}
u(x)=\sum_{n=1}^\infty\langle \pi_n (x)\xi_n, \eta_n\rangle,\quad \xi_n\in E_n,\; \eta_n\in E^*_n,\ x\in G,
\end{equation}
where 
\begin{align*}
(\pi_n,E_n)_{n\in\mathbb{N}}\subseteq\text{Cyc}_p(G),\quad\text{and}\quad \sum_{n=1}^\infty\|\xi_n\| \|\eta_n\|<\infty,
\end{align*}
equipped with the norm
\begin{equation*}
\| u\|=\inf\Big\{\sum_{n=1}^\infty\|\xi_n\| \|\eta_n\| \ :\ u(x)=\sum_{n=1}^\infty\langle \pi_n (x)\xi_n, \eta_n\rangle,\; x\in G\Big\},
\end{equation*}
 is denoted by $B_p(G)$, and is called the $p$-analog of the Fourier-Stieltjes algebra of the locally compact group $G$.
\end{definition}

\begin{remark}\label{REMARKRUNDERUNDE}
\begin{enumerate}
\item\label{REMARKRUNDERUNDE1}
By \cite[Lemma 4.6]{RUNDE2005}, the space $B_p(G)$ can be defined to be the set of all coefficient functions of a $p$-universal representation $(\pi,E)$, and the norm of an element $u\in B_p(G)$ is the infimum of all values $\sum_{n=1}^\infty\|\xi_n\|\|\eta_n\|<\infty$, which such vectors exist in the representation of $u$ as a coefficient function of $(\pi,E)$, i.e. $u(\cdot)=\sum_{n=1}^\infty\langle\pi(\cdot)\xi_n,\eta_n\rangle$.
\item\label{REMARKRUNDERUNDE2} In \cite[Lemma 2.4]{RUNDE2007}, the following identification is shown for an open subgroup $G_0$ of a locally compact group $G$  
\begin{align*}
A_p(G_0)\cong\{f\in A_p(G)\ :\ \text{supp}(f)\subset G_0\},
\end{align*}
and through this fact, one can assume that functions in $A_p(G_0)$ are restriction of functions in $A_p(G)$ to the open subgroup $G_0$.
\end{enumerate}
\end{remark}

\begin{definition}
Let $ (\pi , E)\in \text{Rep}_p(G) $.
\begin{enumerate}
\item For each $f\in L_1(G)$, let $\| f\|_{\pi}:=\|\pi(f)\|_{\mathcal{B}(E)}$, then $\|\cdot\|_{\pi}$ defines an algebra seminorm on $L_1(G)$.
\item By $ PF_{p,\pi}(G) $, we mean the $ p $-pseudofunctions associated with $ (\pi , E) $, which is the closure of
$ \pi( L_1(G)) $ in $ \mathcal{B}(E) $.
\item
If $ (\pi , E) = (\lambda_p,L_p(G))$, we denote $PF_{p,\lambda_p}(G)$ by  $PF_p(G)$.
\item
If $ (\pi , E) $ is $ p $-universal, we denote $ PF_{p,\pi}(G) $ by $ UPF_{p}(G)$, and call it the algebra of
universal $ p $-pseudofunctions.
\end{enumerate}
\end{definition}
\begin{remark}
\begin{enumerate}
\item
For $ p = 2 $, the algebra $ PF_{p}(G) $ is the reduced group $ C^* $-algebra, and $ UPF_{p}(G)$ is the
full group $ C^* $-algebra of $ G $.
\item
If $ (\rho , F)\in \text{Rep}_p(G) $ is such that $ (\pi , E) $ contains every cyclic subrepresentation of $ (\rho , F)$, then $\|\cdot\|_{\rho}\leq\|\cdot\|_{\pi}$ holds. In particular, the definition of $ UPF_{p}(G)$ is independent of a particular $ p $-universal representation.
\item
With $\langle\cdot,\cdot\rangle$ denoting $L_1(G)-L_\infty(G)$ duality, and with $(\pi,E)$ a $p$-universal representation of $G$, we have
\begin{align*}
\| f\|_\pi=\sup\{|\langle f,g\rangle|\ :\ g\in B_p(G),\ \|g\|_{B_p(G)}\leq 1\},\quad f\in L_1(G).
\end{align*}
\end{enumerate} 
\end{remark}

Next lemma states that $ B_p(G) $  is a dual space.
\begin{lemma}{\cite[Lemma 6.5]{RUNDE2005}}\label{RUNDEDUALITY}
Let $ (\pi , E)\in \text{Rep}_{p}(G)$. Then, for each $  \phi\in PF_{p,\pi}(G)^*$, there is a unique $g\in B_p(G) $, with $\|g\|_{B_p(G)}\leq \|\phi\| $ such that
\begin{equation}\label{duality}
\langle\pi (f),\phi\rangle=\int_G f(x)g(x)dx,\qquad f\in L_1(G).
\end{equation}
Moreover, if $ (\pi , E) $ is $p$-universal, we have $\|g\|_{B_p(G)}= \|\phi\| $.
\end{lemma}

\section{Generalized Functorial Properties}\label{SECTIONGENERALIZATIONS}

As the aim of the present paper, here we generalize some properties of the Fourier-Stieltjes algebra, $B(G)$, to the $p$-analog of the Fourier-Stieltjes algebra, $B_p(G)$, for a locally compact group $G$. In the first step, we generalize the idempotent theorem, that was stated in the most possible version as Theorem 1.5 in \cite{RUNDE2007} in Theorem \ref{IDMTHEOREM}. Next, we introduce the $p$-analog of the $\pi$-Fourier and $\pi$-Fourier-Stieltjes spaces given by Arsac in \cite{ARSAC1976} in
Definition \ref{DEFINITIONUNICONV}, and state the expected properties in lemma and proposition afterwards. By such generalization, we then approach to extending result,
Theorem \ref{THEOREMCONCLUSION}, which would be crucial on solving similar problems around $B(G)$ for the $p$-analog case. For instance, in studies on homomorphisms on the $p$-analog of the Fourier-Stieltjes algebra, $B_p(G)$, which can be considered as a $p$-analog of what has been done in \cite{ILIESPRONK2005}, Theorem \ref{THEOREMCONCLUSION} would be essential.

\subsection{Idempotent Theorem}\label{subsectionIT}

In order to being prepared for Theorem \ref{IDMTHEOREM}, which is a generalization of \cite[Theorem 1.5]{RUNDE2007}, we need some elementary definitions and facts, we give them in the following.

\begin{definition}\label{DEFINITIONUNICONV}
\begin{enumerate}
\item
A Banach space $(E, \|\cdot\|)$ is said to be uniformly convex if for
every $0 < \epsilon\leq 2$ there is $\delta > 0$ so that for any two vectors $x$ and $y$ in $E$ with $\|x\|=\|y\| = 1$, the condition $\|x-y\| \geq\epsilon$ implies that $\|\frac{x+y}{2}\|\leq 1-\delta$. Intuitively, the center of a line segment inside the unit ball must lie deep inside
the unit ball unless the segment is short.
\item
A Banach space $E$ is said to be smooth if for each $\xi\in E\backslash\{0\}$ there exists a unique $\eta\in E^*$ such that $\|\eta\|=1$ and $\langle\xi,\eta\rangle =\|\xi\|$.

\end{enumerate}

\end{definition}
\begin{remark}\label{REMUNICONV}
It is worthwhile to note that by Definition \ref{DEFINITIONUNICONV}, every closed
subspace of a uniformly convex Banach space is again a uniformly convex Banach
space.
\end{remark}
Now we state an immensely important theorem about a quotient space which
can be found in \cite{ISTRATESCU1983}.
\begin{theorem}\cite[Theorem 2.4.18]{ISTRATESCU1983}\label{THEOREMQUOTUNI}
Let $E$ be a uniformly convex Banach space and $F$ be a closed linear subspace of $E$. Then the quotient space $E/F$ is uniformly convex Banach space.
\end{theorem}
Now we can conclude the following statement.
\begin{corollary}\label{QSLPSMOOUNICON}
 Every $QSL_p$-space $E$ is uniformly convex and smooth.
 
\begin{proof}
Uniformly convexity of $QSL_p$-space $E$ can be derived from Remark \ref{REMUNICONV}
and Theorem \ref{THEOREMQUOTUNI}. Since $E$ is uniformly convex, by \cite[Lemma 8.4(i) and
Theorem 9.10]{FHHMPZ2001} it is concluded that $E^*$ is smooth, but $E^*$ is a $QSL_{p'}$-space so is
uniformly convex, and then $E^{**}$ is smooth, but $E = E^{**}$ so $E$ is smooth.
\end{proof} 

\end{corollary}

\begin{theorem}\label{IDMTHEOREM}
For a subset $C\subset G$ following statements are equivalent.
\begin{enumerate}
\item\label{idm1}
$C$ is a left open coset,
\item\label{idm2}
$ \chi_C\in B(G) $ with $ \|\chi_C\|_{B(G)}=1 $,
\item\label{idm3}
$\chi_C\neq 0$ is a normalized coefficient function of a representation $(\pi,E)$ where $E$ or $E^*$ is smooth,
\item\label{idm4}
$ \chi_C\in B_p(G) $ with $ \|\chi_C\|_{B_p(G)}=1 $.
\end{enumerate}
\begin{proof}
Equivalency of the first three statements have been proved in \cite[Theorem 1.5]{RUNDE2007}. We demonstrate \eqref{idm2}$\Rightarrow$\eqref{idm4}$\Rightarrow$\eqref{idm3}. Let \eqref{idm2} hold. Then from the fact that $B(G)\subset B_p(G)$ and this embedding is a contraction, we have $\chi_C\in B_p(G)$ with $ \|\chi_C\|_{B_p(G)}\leq 1$, which by inequality $\|\cdot\|_{C_b(G)}\leq \|\cdot\|_{B_p(G)}$, we have $ \|\chi_C\|_{B_p(G)}=1$ which shows \eqref{idm2} implies \eqref{idm4}.\\
Now let $\chi_C\in B_p(G)$ with $ \|\chi_C\|_{B_p(G)}= 1$. So, by Definition \ref{definitionofBp}, the function $\chi_C$ is a normalized coefficient function of an isometric group representation on a $QSL_{p}$-space, which is smooth by Corollary \ref{QSLPSMOOUNICON} that is \eqref{idm3}.
\end{proof}

\end{theorem}

\begin{corollary}\label{COROIDEMCOSETRING}
Let $G$ be a locally compact group and $Y\in\Omega_0(G)$, then we have $\chi_Y\in B_p(G)$. Moreover, we have
\begin{align}\label{FFF}
1\leq\|\chi_Y\|_{B_p(G)}\leq 2^{m_Y},\quad\text{with}\quad m_Y=\inf\{m\in\mathbb{N}\; :\ Y=Y_0\backslash\cup_{i=1}^mY_i\},
\end{align}
where for $i=0,1,\ldots,m$ sets $Y_i$, are as \eqref{OPENCOSETRINGS}.
\begin{proof}
Since $Y\in\Omega_0(G)$, then by \eqref{OPENCOSETRINGS}, there exist open coset $Y_0$ and open subcosets $Y_i\subset Y_0$, for $i=1,\ldots , m$ and $m\in\mathbb{N}$ such that $Y=Y_0\backslash\cup_{i=1}^mY_i$. By Theorem \ref{IDMTHEOREM}-\eqref{idm4}, we have $\chi_{Y_i}\in B_p(G)$, with $\|\chi_{Y_i}\|_{B_p(G)}=1$, for $i=0,1,\ldots ,m$. On the other hand, since 
\begin{align}\label{EQCHI}
\chi_Y=\chi_{Y_0}-\big(\sum_{i=1}^m\chi_{Y_i}-\sum_{i,j=1}\chi_{Y_i\cap Y_j}+\sum_{i,j,k=1}\chi_{Y_i\cap Y_j\cap Y_k}+\cdots+(-1)^{m+1}\chi_{Y_1\cap Y_2\cap\ldots\cap Y_m}\big),
\end{align}
then we have $\|\chi_Y\|_{B_p(G)}\leq 2^{m_Y}$, and by taking infimum on all possible decomposition of $Y$ as \eqref{OPENCOSETRINGS} relation \eqref{FFF} holds.
\end{proof}

\end{corollary}

\subsection{$p$-Analog of the $\pi$-Fourier spaces}\label{subsectionAFS}

In the sequel, we will give some extensions of results in \cite{ARSAC1976}. For a unitary representation $(\pi,\mathcal{H}_\pi)$ with Hilbert space $\mathcal{H}_\pi$, the $\pi$-Fourier space  has been defined to be closed linear span of the set of the coefficient functions of the representation $(\pi,\mathcal{H}_\pi)$, and is denoted by $A_\pi$, with the norm in usual way. Moreover, $\pi$-Fourier-Stieltjes algebra, $B_\pi$, for such representation is defined to be $w^*$-closure of $A_\pi$. Additionally, if we let $C^*_\pi(G)$ be the  $C^*$-algebra associated with $\pi$, we have $B_\pi=C^*_\pi(G)^*$. Here we introduce $p$-generalization of these results.

\begin{definition}
For a representation $(\pi , E)\in\text{Rep}_p(G)$, we define the \textit{$p$-analog of the $\pi$-Fourier space}, $A_{p,\pi}$, to be closed linear span of the collection of the coefficient functions of representation $(\pi, E)$, i.e. functions of the form
\begin{align*}
u(x)=\sum_n\langle\pi(x)\xi_n,\eta_n\rangle,\quad x\in G, (\xi_n)_{n\in\mathbb{N}}\subseteq E, \; (\eta_n)_{n\in\mathbb{N}}\subseteq E^*
\end{align*}
 equipped with the norm 
\begin{align*}
\| u\|_{A_{p,\pi}}=\inf\Big\{\sum_{n=1}^\infty\|t_n\|\|s_n\|\ :\  u(x)=\sum_n\langle\pi(x)t_n,s_n\rangle,\ x\in G\Big\},
\end{align*}
and evidently, infimum is taken over all possible equivalent representative of $u$ so that the value is convergent.
\end{definition}

\begin{remark}\label{REMARKPSIMAP}
\begin{enumerate}
\item\label{REMARKPSIMAP1}
For $(\pi , E)\in\text{Rep}_p(G)$ consider the map $\Psi_{p,\pi}: E^*\widehat{\otimes}E\rightarrow C_b(G)$, defined via
\begin{align*}
\Psi_{p,\pi}\bigg(\sum_{n}\eta_n\otimes\xi_n\bigg)=\sum_{n}\langle\pi(x)\xi_n,\eta_n\rangle,\ x\in G.
\end{align*}
This map is onto to its range, which is $A_{p,\pi}$, and  so we can identify it with the Banach space $E^*\widehat{\otimes}E/\ker{\Psi_{p,\pi}}$, and the norm on $A_{p,\pi}$ is the quotient norm i.e.,
\begin{align*}
\Big\| \sum_{n}\xi_n\otimes\eta_n+\ker\Psi_{p,\pi}\Big\|&=\inf\Big\{\sum_{n}\|t_n \| \| s_n \| \ : \sum_{n}\langle\pi(\cdot)t_n,s_n\rangle=\sum_{n}\langle\pi(\cdot)\xi_n,\eta_n\rangle \Big\}\\
&=\Big\|\sum_{n}\langle\pi(\cdot)\xi_n,\eta_n\rangle \Big\|_{A_{p,\pi}}.
\end{align*}
So, one can identify $A_{p,\pi}$ with the quotient space $E^*\widehat{\otimes}E/\ker\Psi_{p,\pi}$.
\item
Since we have $A_{p,\pi}\cong E^*\widehat{\otimes}E/\ker\Psi_{p,\pi}$, then the space $A_{p,\pi}$ is a Banach space.

\end{enumerate}
\end{remark}

In the next proposition we give an equivalent formula of computing norm on the space $A_{p,\pi}$, for a representation $(\pi,E)\in \text{Rep}_p(G)$. For this aim, we denote the set of cyclic subrepresentations of a representation $(\pi,E)$ by $\text{Cyc}_{p,\pi}(G)$.

\begin{proposition}\label{LEMMANORMAPP}

Let $(\pi , E)\in\text{Rep}_p(G)$ and  $u\in A_{p,\pi}$. Then we have
\begin{align}\label{eq11}
\|u\|_{A_{p,\pi}}=\inf\Big\{\sum_{n}\|t_n \| \| s_n \| \ : \ u(x)=\sum_{n}\langle\rho_n(x)t_n,s_n\rangle,\ x\in G\Big\},
\end{align}
where the infimum is taken on all representations of $u$, in which $\Big((\rho_n, F_n)\Big)_{n\in\mathbb{N}}\subseteq\text{Cyc}_{p,\pi}(G)$ with $(t_n)_{n\in\mathbb{N}}\subseteq F_n $ and $(s_n)_{n\in\mathbb{N}}\subseteq F^*_n$.
\begin{proof}
Let us denote the infimum in \eqref{eq11} by $C$. Assume that for $x\in G$, we have $u(x)=\sum_{n}\langle\pi(x)\xi_n,\eta_n\rangle$ with $\sum_{n}\|\xi_n\|\|\eta_n\|<\infty$. For each $n\in\mathbb{N}$, if we put
\begin{align*}
&F_n=\overline{\pi(L_1(G))\xi_n{}{}}^{\|\cdot\|_{E}},\quad \rho_n :G\rightarrow\mathcal{B}(F_n),\quad \rho_n(x)=\pi(x)|_{F_n},\ x\in G,
\end{align*}
and $t_n=\xi_n,\ s_n=\eta_n|_{F_n}$, then we have
\begin{align*}
&\Big((\rho_n, F_n)\Big)_{n\in\mathbb{N}}\subseteq\text{Cyc}_{p,\pi}(G),\quad u(x)=\sum_{n}\langle\rho_n(x)t_n,s_n\rangle,
\end{align*}
with $ C\leq\sum_{n=1}^\infty\|t_n\|\|s_n\|\leq \sum_{n}\|\xi_n\|\|\eta_n\|$.
Since $(\xi)_{n\in\mathbb{N}}\subset E$ and $(\eta)_{n\in\mathbb{N}}\subset E^*$ are arbitrary in the representing of $u$, we have $C\leq\| u\|_{A_{p,\pi}}$.\\
For the inverse inequality, let $\epsilon >0$ is given. Then there exist $((\rho_n, F_n))_{n\in\mathbb{N}}\subseteq\text{Cyc}_{p,\pi}(G)$, $(t_n)_{n\in\mathbb{N}}\subseteq F_n$, $(s_n)_{n\in\mathbb{N}}\subseteq F_n^*$, such that for each $n\in\mathbb{N}$, we have $(\rho_n, F_n)\subset (\pi,E)$ and
\begin{equation*}
\sum_{n}\|t_n\|\|s_n\| < C+\epsilon,\qquad u(x)=\sum_{n}\langle\rho_n(x)t_n,s_n\rangle,\quad x\in G.
\end{equation*}
Now for each $n\in\mathbb{N}$, by applying Hahn-Banach theorem extend each $s_n\in F_n^*$ to the $\eta_n\in E^*$ such that $\|\eta_n\|=\|s_n\|$. Therefore,
\begin{equation*}
\| u\|_{A_{p,\pi}}\leq \sum_{n}\|t_n\|\|\eta_n\| =\sum_{n}\|t_n\|\|s_n\| < C+\epsilon,
\end{equation*}
and it means $\| u\|_{A_{p,\pi}}\leq C$.
\end{proof}
\end{proposition}
In the next attempt, we will generalize some previously existed result, and for this aim, we need to be precise about notions there. To clarify everything, we bring some definitions and their requirement in the sequel. To do this, we briefly state some facts about ultrafilters on Banach spaces. The main reference here is \cite{HEINRICH1986}, and the more applicable one is \cite{DAWS2004}.\\
Let $(E_i)_{i\in \mathbb{I}}$ be a collection of Banach spaces for an index set $\mathbb{I}$ . Consider the Banach space $L_\infty(\mathbb{I} , E_i)$ of elements $(\xi_i)_{i\in \mathbb{I}}$ equipped with the pointwise operations and the supremum norm $\|(\xi_i)_i\|= \sup_{i\in \mathbb{I}}\|\xi_i\|_{E_i}$. Now, let $\mathcal{U}$ be an ultrafilter on $\mathbb{I}$ , and let $N_\mathcal{U}$ be as following
\begin{align*}
N_\mathcal{U} = \Big\{(\xi_i)_{i\in \mathbb{I}}\  :\  \lim_\mathcal{U}\|\xi_i\|=0\Big\},
\end{align*}
where by $\lim_\mathcal{U}\|\xi_i\|$ we mean the limit of $(\xi_i)_{i\in \mathbb{I}}$ along the ultrafilter $\mathcal{U}$,
that exists due to the fact that the values $\|\xi_i\|$ belong to the compact interval $[0, M]$, where $\sup_{i\in \mathbb{I}} \|\xi_i\|= M <\infty$. Note that $N_\mathcal{U}$
is a closed subspace, therefore, the completion of the quotient space $L_\infty(\mathbb{I} , E_i)/N_\mathcal{U}$ is a Banach space. This Banach space is denoted by $(E_i)_\mathcal{U}$ and it is called the ultraproduct of $(E_i)_{i\in \mathbb{I}}$ with respect to the ultrafilter
$\mathcal{U}$. Besides, the quotient norm of an element $(\xi_i)_\mathcal{U}\in (E_i)_{\mathcal{U}}$ coincides with
the $\lim_\mathcal{U}\|\xi_i\|$.\\
In general, we have $(E^*_i)_\mathcal{U}\subseteqq (E_i)_\mathcal{U}^*$, isometrically. In addition, when for each $i\in I$ we have $E_i = E$, then the ultraproduct space is called the ultrapower of the Banach space $E$ and is denoted by $(E)_\mathcal{U}$ .  Moreover, if $E$ is a super-reflexive space, then we have $(E_i)_\mathcal{U}^*=(E_i^*)_\mathcal{U}$ , and this is a well-known fact that an ultrapower of a $QSL_p$-space $E$ is again a $QSL_p$-space.\\
For a Banach space $E$, the natural embedding $J : E \rightarrow (E)_\mathcal{U}$ , defined via $J(\xi) =(\xi_i)_\mathcal{U}$ where $\xi_i\equiv\xi$, is an isometric one. Furthermore, for the case that $E$ is super-reflexive, if $(\eta_i)_\mathcal{U}\in (E^*)_\mathcal{U}=(E)^*_\mathcal{U}$ then we have 
\begin{align*}
\langle J(\xi),(\eta_i)_\mathcal{U}\rangle=\lim_\mathcal{U}\langle\xi,\eta_i\rangle.
\end{align*}
Additionally, there is a well-known canonical isometric map $\kappa_E : E \rightarrow E^{**}$ that is defined by $\langle\eta,\kappa_E(\xi)\rangle =\langle\xi,\eta\rangle$, and it is a surjection, if and only if $E$ is reflexive.
For a Banach space $E$, and an ultrafilter $\mathcal{U}$, since the unit ball of $E^{**}$ is compact, then the following map is well-defined and contractive
\begin{align*}
\mathcal{J}:(E)_\mathcal{U}\rightarrow E,\quad \mathcal{J}((\xi_i)_\mathcal{U})=w^*-\lim_\mathcal{U}\kappa_E(\xi_i),
\end{align*}
and for every $\eta\in E$ we have $\langle\eta,\mathcal{J}((\xi_i)_\mathcal{U})\rangle=\lim_\mathcal{U}\langle\xi_i,\eta\rangle$.\\
According to \cite[Proposition 2]{DAWS2004}, for a Banach space $E$ there exists an ultrafilter $\mathcal{U}$ and an isometric map, $\overline{J}:E^{**}\rightarrow (E)_\mathcal{U}$ such that $\overline{J}|_E=J$, $\mathcal{J}\circ\overline{J}=\mathrm{id}_{E^{**}}$, and $\overline{J}\circ\mathcal{J}$ is a norm one projection from $(E)_\mathcal{U}$ onto $\overline{J}(E^{**})$.\\
For a Banach space E, and a complex number $p\in (1,\infty)$, the Banach space $L_p(E) = L_p(\mathbb{N} , E)$, is as following
\begin{align*}
L_p(E)=\Big\{(\xi_n)_n\ :\ \|(\xi_n)_n\|=\Big(\sum_n\|\xi_n\|^p\Big)^\frac{1}{p}<\infty\Big\},
\end{align*}
which by \cite[Proposition 4]{DAWS2004} it is super-reflexive whenever $E$ is so. Additionally, in the case that $E$ is a $QSL_p$-space, then $L_p(E)$ is again a $QSL_p$-space.\\
Now, for a representation $(\pi,E)\in\text{Rep}_p(G)$, we introduce two representations $(\pi^\infty, L_p(E))$ and $(\pi_\mathcal{U}, (E)_\mathcal{U})$, as following:
\begin{align*}
&\pi^\infty:G\rightarrow\mathcal{B}(L_p(E)),\quad \pi^\infty(x)((\xi_n)_n)=(\pi(x)\xi_n)_n,\quad x\in G,\; (\xi_n)_n\in L_p(E),\\
&\pi_\mathcal{U}:G\rightarrow\mathcal{B}((E)_\mathcal{U}),\quad \pi_\mathcal{U}(x)((\xi_i)_\mathcal{U})=(\pi(x)\xi_i)_\mathcal{U},\quad x\in G,\; (\xi_i)_\mathcal{U}\in (E)_\mathcal{U}.
\end{align*}

\begin{lemma}\label{LEMMANEWREP}
For a representation $(\pi,E)\in\text{Rep}_p(G)$, consider the aforementioned representations $(\pi^\infty, L_p(E))$ and $(\pi_\mathcal{U}, (E)_\mathcal{U})$. Then the following statements hold.
\begin{enumerate}
\item\label{LEMMANEWREP1}
The ranges of both representations are the subspaces of invertible and isometric operators on their associated $QSL_p$-spaces, and therefore, they belong to $\text{Rep}_p(G)$.
\item\label{LEMMANEWREP2}
These representations are related to $(\pi,E)$ as representations of the group algebra $L_1(G)$ as the same as they are related as representations of the group $G$. Precisely, for an element $f\in L_1(G)$, we have
\begin{align*}
\pi^\infty(f)=(\pi(f))^\infty,\qquad \pi_\mathcal{U}(f)=(\pi(f))_\mathcal{U}.
\end{align*}
\item\label{LEMMANEWREP3}
$(\pi^\infty, L_p(E))$ and $(\pi_\mathcal{U}, (E)_\mathcal{U})$ are essential representations of $L_1(G)$.
\end{enumerate}
\begin{proof}
The two first parts hold naturally. For the third part we briefly assure the  reader about our claim.
\item{\textit{Case one:} $(\pi^\infty,L_p(E))$.} Let $(\xi_n)_n\in L_p(E)$. Since $(\pi, E)$ is an essential representation, then for every given $\epsilon>0$, for each $n\in\mathbb{N}$, there exist  $f_n\in L_1(G)$ and $t_n\in E$ such that
\begin{align*}
\|\pi(f_n)t_n-\xi_n\|<\frac{\epsilon}{nM},\quad M=\Big(\sum_n\frac{1}{n^p}\Big)^\frac{1}{p}<\infty\ (\text{since}\; p>1),
\end{align*}
then we have
\begin{align*}
\|(\pi(f_n)t_n)_n-(\xi_n)_n\|<\epsilon.
\end{align*}
So, the arbitrary element $(\xi_n)_n$ is approximated by the element $(\pi(f_n)t_n)_n$, that lives in the range of $\pi^\infty$, as a representation of $L_1(G)$.
\item{\textit{Case two:} $(\pi_\mathcal{U},(E)_\mathcal{U})$.} Let $(\xi_i)_\mathcal{U}\in (E)_\mathcal{U}$, and consider a representative $(\xi_i)_{i\in\mathbb{I}}\in L_\infty(\mathbb{I},E)$ of it. By the same argument as above, for every given $\epsilon>0$, for each $i\in\mathbb{I}$, there exist $f_i \in L_1(G)$, and $t_i\in E$ such that
\begin{align*}
\|\pi(f_i)t_i-\xi_i\|<\frac{\epsilon}{2},
\end{align*}
and consequently,
\begin{align*}
\|(\pi(f_i)t_i)_i-(\xi_i)_i\|\leq\frac{\epsilon}{2}<\epsilon.
\end{align*}
On the other hand, since for an element $(r_i)_i\in L_\infty(\mathbb{I},E)$ we have $\|(r_i)_\mathcal{U}\|\leq\|(r_i)_i\|$ it is obtained that
\begin{align*}
\|(\pi(f_i)t_i)_\mathcal{U}-(\xi_i)_\mathcal{U}\|<\epsilon.
\end{align*}
Throughout the fact that the element $(\pi(f_i)t_i)_\mathcal{U}$ belongs to the range of $\pi_\mathcal{U}$ as representation of $L_1(G)$, we are done.
\end{proof}
\end{lemma}

\begin{corollary}
For a representation $(\pi,E)\in\text{Rep}_p(G)$, the representation $(\pi^\infty_\mathcal{U},(L_p(E))_\mathcal{U}$, defined in the obvious way, is an essential representation and belongs to $\text{Rep}_p(G)$.
\end{corollary}

Next proposition is somewhat a restatement of the Lemma \ref{RUNDEDUALITY} which Lemma 6.5 in \cite{RUNDE2005}. This restatement is beneficial due to detailed proof and 

\begin{proposition}\label{PROPULTRAAP}
Let $(\pi,E)\in\text{Rep}_p(G)$. Then
\begin{enumerate}
\item
there exists a free ultrafilter $\mathcal{U}$, such that the canonical representation of $PF_{p,\pi}(G)$ on $F = (L_p(E))_\mathcal{U}$ is weak-weak$^*$ continuous, essential and
isometric,
\item the identification $PF_{p,\pi}(G)^*=\overline{A_{p,\pi^\infty_\mathcal{U}}{}}^{w^*}= A_{p,\pi^\infty_\mathcal{U}}$
holds.
\end{enumerate}
\begin{proof}
Here, we sometimes use $F$ instead of $(L_p(E))_\mathcal{U}$ , for ease of notation, and sometimes do not use to highlight the associated space, and actions.\\
To prove part one, by the proof of \cite[Proposition 5]{DAWS2004}, there exists an ultraﬁlter $\mathcal{U}$ on an indexing set $\mathbb{I}$, such that by considering the above mentioned map $\bar{J}$ for $E^*\widehat{\otimes}E$ and using the fact that $(E^*\widehat{\otimes}E)^*=\mathcal{B}(E)$, we have
\begin{align*}
J : \mathcal{B}(E)^*\rightarrow (E^*\widehat{\otimes}E)_\mathcal{U}
\end{align*}
isometrically. Through \cite[Theorem 1 or Proposition 5]{DAWS2004}, for the obtained ultrafilter above, the map $P: F^*\widehat{\otimes}F\rightarrow\mathcal{B}(E)^*$ defined for $t =((t_{i,n})_n)_\mathcal{U}\in F$ and $s =((s_{i,n})_n)_\mathcal{U}\in F^*$ via
\begin{align*}
\langle T, P(s\otimes t)\rangle=\lim_\mathcal{U}\sum_n\langle T(t_{i,n}),s_{i,n}\rangle,\quad T\in\mathcal{B}(E).
\end{align*}
is a linear metric surjection. Therefore, the embedding $P^*:\mathcal{B}(E)^{**}\rightarrow\mathcal{B}(F)$ is an isometric homomorphism. So, the canonical representation of $PF_{p,\pi}(G)\subseteq \mathcal{B}(E)$ on $F = (L_p(E))_\mathcal{U}$ that is $P^*|_{PF_{p,\pi}(G)} = P^*_r$, is weak-weak$^*$ and isometric. Precisely, the following map satisfies mentioned properties
\begin{align*}
P^*_r:PF_{p,\pi}(G)\rightarrow \mathcal{B}(F),\quad P^*_r(\pi(f))=\pi^\infty_\mathcal{U}(f),\ f\in L_1(G).
\end{align*}
In fact, as an application of \ref{LEMMANEWREP}-\eqref{LEMMANEWREP3}, through the following diagram,
\begin{displaymath}
\xymatrix{ PF_{p,\pi}(G)\ar[r]^{P^*_r} & \mathcal{B}(L_p(E))_\mathcal{U}) \\
L_1(G)\ar[u]^\pi \ar[ur]_{\pi^\infty_\mathcal{U}} &  }
\end{displaymath}
is an essential representation of $PF_{p,\pi}(G)$. In detail, the map $\pi$ is contractive with dense range while $\pi^\infty_\mathcal{U}$ is an essential representation as it is described at the end of previous lemma. Since we have $P^*_r\circ\pi=\pi^\infty_\mathcal{U}$, then our claim is true and $P^*_r$ is an isometric, weak-weak$^*$ continuous, and essential representation of $PF_{p,\pi}(G)$. Subsequently, we have
\begin{align*}
PF_{p,\pi}(G)=PF_{p,\pi^\infty_\mathcal{U}}(G)
\end{align*}
For the second part, since $(P^*_r)^*_r$, the restriction of the conjugate map $(P^*_r)^*$ to the subspace $F^*\widehat{\otimes}F$ is a quotient map onto $PF_{p,\pi}(G)^*$, then we have
\begin{align*}
(P^*_r)^*_r: F^*\widehat{\otimes}F\rightarrow PF_{p,\pi}(G)^*.
\end{align*}
It is obtained that
\begin{align*}
PF_{p,\pi}(G)^*=(L_{p'}(E^*))_\mathcal{U}\widehat{\otimes}(L_{p}(E))_\mathcal{U}/\ker(P^*_r)^*_r.
\end{align*}
For an element $\phi\in PF_{p,\pi}(G)^*$, there exists a unique $\tau\in F^*\widehat{\otimes}F/\ker(P^*_r)^*_r$ such that for a given $\epsilon>0$ there exist $(\xi_k)_k\subseteq F$ and $(\eta_k)_k\subseteq F^*$ with $\tau=\sum_k\xi_k\otimes\eta_k$, and
\begin{align*}
\|\phi\|\leq\sum_k\|\xi_k\|\|\eta_k\|<\|\phi\|+\epsilon.
\end{align*}
Additionally, for every $f\in L_1(G)$, we have
\begin{align}\label{eq22}
\langle\pi(f),\phi\rangle=\sum_k\langle P^*_r\circ\pi(f)\xi_k,\eta_k\rangle=\sum_k\langle\pi^\infty_\mathcal{U}(f)\xi_k,\eta_k\rangle=\langle\pi^\infty_\mathcal{U}(f),u\rangle
\end{align}
where
\begin{align*}
u(x)=\sum_k\langle\pi^\infty_\mathcal{U}(x)\xi_k,\eta_k\rangle\in A_{p,\pi^\infty_\mathcal{U}},\qquad x\in G,
\end{align*}
and $\langle\pi^\infty_\mathcal{U},u\rangle$ means the $L_1-L_\infty$ duality between $f$ and $u$, as it is described in Lemma \ref{RUNDEDUALITY}. Now, consider the map $\Psi_{p,\pi^\infty_\mathcal{U}}$, as Remark \ref{REMARKPSIMAP}-\eqref{REMARKPSIMAP1}, associated to the representation $(\pi^\infty_\mathcal{U},F)$. We have
\begin{align*}
A_{p,\pi^\infty_\mathcal{U}}=(L_{p'}(E^*))_\mathcal{U}\widehat{\otimes}(L_p(E))_\mathcal{U}/\ker\Psi_{\pi,\pi^\infty_\mathcal{U}}.
\end{align*}
Since, $PF_{p,\pi}(G)=PF_{p,\pi^\infty_\mathcal{U}}(G)$, then the relation \eqref{eq22} reveals that kernels of the maps $(P^*_r)^*_r$ and $\Psi_{p,\pi^\infty_\mathcal{U}}$ coincide and we are done.
\end{proof}

\end{proposition}

\begin{proposition}
Let $(\pi,E)\in\text{Rep}_p(G)$. Then we have the following
identification
\begin{align*}
A_{p,\pi}=A_{p,\pi^\infty}.
\end{align*}
\begin{proof}
Let $u\in A_{p,\pi^\infty}$ , and $\epsilon>0$ is given. There exist the sequence of vectors $((\xi_{n,m})_n)_m\subseteq L_p(E)$ and $((\eta_{n,m})_n)_m\subseteq L_p(E)^*=L_{p'}(E^*)$ such that
\begin{align}\label{eq33}
u(x)=\sum_m\langle\pi^\infty(x)(\xi_{n,m})_n,(\eta_{n,m})_n\rangle=\sum_{m,n}\langle\pi(x)\xi_{n,m},\eta_{n,m}\rangle
\end{align}
and
\begin{align}\label{eq44}
\|u\|_{\pi^\infty}+\epsilon>\sum_m\|(\xi_{n,m})_n\|\|(\eta_{n,m})_n\|\geq\sum_{n,m}\|\xi_{n,m}\|\|\eta_{n,m}\|
\end{align}
In the last inequality, we utilized the H\"older inequality of positive numbers.
From \eqref{eq33} it is evident that $u\in A_{p,\pi}$ and \eqref{eq44} shows that $\|u\|_{p,\pi}\leq\|u\|_{p,\pi^\infty}$,
which means that, $A_{p,\pi^\infty}\subseteq A_{p,\pi}$, contractively. We shall show the inverse inclusion holds contractively. To do so, let $u\in A_{p,\pi}$, and for a given $\epsilon>0$ let vectors $(t_n)_n\subseteq E$ and $(s_n)_n\subseteq E^*$ be such that
\begin{align*}
u(x)=\sum_n\langle\pi(x)t_n,s_n\rangle,\qquad x\in G,
\end{align*}
and
\begin{align*}
\|u\|+\epsilon>\sum_{n}\|t_n\|\|s_n\|,
\end{align*}
Now, if we put
\begin{align*}
\xi_n=\|t_n\|^{-1+\frac{1}{p}}\|s_n\|^\frac{1}{p}t_n,\qquad \eta_n=\|s_n\|^{-1+\frac{1}{p'}}\|t_n\|^\frac{1}{p'}s_n
\end{align*}
then we have $(\xi_n)_n\in L_p(E)$ and $(\eta_n)_n\in L_{p'}(E^*)$. Moreover,
\begin{align*}
&\Big(\sum_n\|\xi_n\|^p\Big)^\frac{1}{p}=\Big(\sum_n\|t_n\|\|s_n\|\Big)^\frac{1}{p}<\Big(\|u\|_{\pi}+\epsilon\Big)^\frac{1}{p}\\
&\Big(\sum_n\|\eta_n\|^{p'}\Big)^\frac{1}{p'}=\Big(\sum_n\|t_n\|\|s_n\|\Big)^\frac{1}{p'}<\Big(\|u\|_{\pi}+\epsilon\Big)^\frac{1}{p'}
\end{align*}
and
\begin{align*}
u(x)=\langle\pi^\infty(x)(\xi_n)_n,(\eta_n)_n\rangle=\sum_n\langle\pi(x)\xi_n,\eta_n\rangle=\sum_n\langle\pi(x)t_n,s_n\rangle,\quad x\in G.
\end{align*}
Therefore,
\begin{align*}
u\in A_{p,\pi^\infty},\quad \text{and}\quad \|u\|_{\pi^\infty}\leq\|(\xi_n)_n\|\|(\eta_n)_n\|<\|u\|_\pi+\epsilon.
\end{align*}
\end{proof}
\end{proposition}
\begin{corollary}
For a representation $(\pi,E)\in\text{rep}_p(G)$ we have the following identification
\begin{align*}
PF_{p,\pi}(G)^*=A_{p,\pi_\mathcal{U}}.
\end{align*}
\begin{proof}
It is a straightforward.
\end{proof}
\end{corollary}

\begin{remark}\label{REMARKOFDUALOFPFP}
\begin{enumerate}
\item\label{REMARKOFDUALOFPFP1}
In the light of previous proposition, due to the fact that $(\pi^\infty_\mathcal{U},F)$ is weak-weak$^*$ continuous, essential an isometric representation of $PF_{p,\pi}(G)$, then we have
\begin{align*}
PF_{p,\pi}(G)=PF_{p,\pi^\infty_\mathcal{U}}(G).
\end{align*}
\item\label{REMARKOFDUALOFPFP2}
In the case that the representation $(\pi,E)$ is a $p$-universal representation, then since $(\pi^\infty_\mathcal{U},(l_p(E))_\mathcal{U})$ is also a $p$-universal representation, our notation coincides with Runde's one in \cite{RUNDE2005}.
\item\label{REMARKOFDUALOFPFP3}
We follow \cite{ARSAC1976} in notation, and denote $A_{p,\pi_\mathcal{U}}$ by $B_{p,\pi}$, and we call it \textit{$p$-analog of the $\pi$-Fourier-Stieltjes algebra}, which by Proposition \ref{PROPULTRAAP} is the dual space of the space of $p$-pseudofunctions  associated with $(\pi,E)\in\text{Rep}_p(G)$, i.e. the dual space of $PF_{p,\pi}(G)$ through following duality
\begin{align*}
\langle \pi(f),u\rangle=\int_G u(x)f(x)dx,\quad f\in L_1(G),\ u\in B_{p,\pi},
\end{align*}
and as we expect that, we have
\begin{align*}
&\| u\|=\sup_{\|f\|_\pi\leq 1}|\langle \pi(f), u\rangle|=\sup_{\|f\|_\pi\leq 1}|\int_G u(x)f(x)dx|,\quad u\in B_{p,\pi},\\
&\| f\|_\pi=\sup_{\|u\|\leq 1}|\langle \pi(f), u\rangle|=\sup_{\|u\|\leq 1}|\int_G u(x)f(x)dx|,\quad f\in L_1(G).\\
\end{align*}
So, we have set $PF_{p,\pi}(G)^*=A_{p,\pi_\mathcal{U}}=B_{p,\pi}$ , and in the case that $(\pi,E)=(\lambda_{p,G},L_p(G))$ we usually use the symbol $PF_{p}(G)^*$.

\item\label{REMARKOFDUALOFPFP4}
It is obvious that $B_{p,\pi}\subseteq B_p(G)$ is a contractive inclusion for every $(\pi , E)\in\text{Rep}_p(G)$, and if $(\pi,E)$ is a $p$-universal representation it will become an isometric isomorphism.
\item\label{REMARKOFDUALOFPFP5}
It is valuable to note that the ultra filter $\mathcal{U}$ is the one for which the embedding $\mathcal{B}(E)^{*}\subseteq (E^*\widehat{\otimes}E)_\mathcal{U}$ is isometric, so the space $A_{p,\pi^\infty_\mathcal{U}}$ is determined. Furthermore, if $\mathcal{V}$ is another free ultrafilter that makes the similar embedding $\mathcal{B}(E)^{*}\subseteq (E^*\widehat{\otimes}E)_\mathcal{V}$ into an isometry, then we have
\begin{align*}
A_{p,\pi_\mathcal{U}}=PF_{p,\pi}(G)^*= A_{p,\pi_\mathcal{V}}.
\end{align*}
So, our definition is independent of choosing suitable free ultrafilter, therefore, it is well-defined.
\item\label{REMARKOFDUALOFPFP6}
For a locally compact group $G$, we have the following contractive inclusions
\begin{align*}
PF_{p}(G)^*=B_{p,\lambda_p}\subset B_p(G)\subset \mathcal{M}(A_p(G)).
\end{align*}
All inclusions will become equalities in the case that $G$ is amenable (See \cite[Theorem 6.6 and Theorem 6.7]{RUNDE2005}).
\end{enumerate}
\end{remark}

\subsection{Extension Theorem}\label{subsectionET}

In the following, we study some functorial properties of the $p$-analog of the Fourier-Stieltjes algebras. One of the earliest questions about such algebras is when an extension of a function defined on a subgroup belongs to the $p$-analog of the Fourier-Stieltjes algebra on the larger group . For this aim, we deal with following notation. Let $G_0\subset G$, be any subset, and $u:G_0\rightarrow \mathbb{C}$ be a function. By $u^\circ$ we mean
\begin{align*}
u^\circ=\left\{
\begin{array}{lll}
u&&\text{on}\; G_0\\
0&&\text{o.w.}
\end{array}\right. .
\end{align*}
Next lemma is going to express the relation between representation of an open subgroup $G_0$ with the one of the initial group $G$.

\begin{lemma}\label{LEMMARESTRICTIONMAP}
Let $(\pi,E)\in\text{Rep}_p(G)$. Then the restriction of $\pi$ to the open subgroup $G_0$, which is denoted by $(\pi_{G_0},E)$ belongs to $\text{Rep}_p(G_0)$. Moreover, for each $f\in L_1(G_0)$ and each $g\in L_1(G)$, we have the following relations
\begin{align}\label{w0}
\pi_{G_0}(f)=\pi(f^\circ),\quad\text{and}\quad \pi_{G_0}(g|_{G_0})=\pi(g\chi_{G_0}).
\end{align}
\begin{proof}
It is evident that $(\pi_{G_0},E)\in\text{Rep}_p(G_0)$. For the second part, simple calculations below reveal that our claim is true. For $\xi\in E$ and $\eta\in E^*$, if $f\in L_1(G_0)$ and $g\in L_1(G)$, then we have
\begin{align*}
\langle\pi_{G_0}(f)\xi,\eta\rangle &=\int_{G_0}f(x)\langle\pi_{G_0}(x)\xi,\eta\rangle dx\\
&=\int_{G_0}f(x)\langle\pi(x)\xi,\eta\rangle dx\\
&=\int_{G}f^\circ(x)\langle\pi(x)\xi,\eta\rangle dx\\
&=\langle\pi(f^\circ)\xi,\eta\rangle,
\end{align*}
and
\begin{align*}
\langle\pi_{G_0}(g|_{G_0})\xi,\eta\rangle &=\int_{G_0}g|_{G_0}(x)\langle\pi_{G_0}(x)\xi,\eta\rangle dx\\
&=\int_{G_0}g|_{G_0}(x)\langle\pi(x)\xi,\eta\rangle dx\\
&=\int_{G}g(x)\chi_{G_0}(x)\langle\pi(x)\xi,\eta\rangle dx\\
&=\langle\pi(g\chi_{G_0})\xi,\eta\rangle.
\end{align*}
So, we have
\begin{align}
&\langle\pi_{G_0}(f)\xi,\eta\rangle =\langle\pi(f^\circ)\xi,\eta\rangle,\label{w1}\\
&\langle\pi_{G_0}(g|_{G_0})\xi,\eta\rangle =\langle\pi(g\chi_{G_0})\xi,\eta\rangle.\label{w2}
\end{align}
and since \eqref{w1} and \eqref{w2} hold for every $\xi\in E$ and $\eta\in E^*$, then the relations in \eqref{w0} are obtained.
\end{proof}
\end{lemma}

\begin{proposition}\label{PROPEXREC}
Let $G$ be a locally compact group and $G_0$ be its open
subgroup, and let $(\pi,E)\in \text{Rep}_p(G)$. Then the following statements hold.
\begin{enumerate}
\item\label{PROPEXREC1}
The map $S_{\pi_{G_0}}:PF_{p,\pi_{G_0}}(G_0)\rightarrow PF_{p,\pi}(G)$ defined via $S_{\pi_{G_0}}(\pi_{G_0}(f))=\pi(f\chi_{G_0})$, for $f \in L_1(G_0)$ is an isometric homomorphism. In fact, we have
the following isometric identification
\begin{align*}
PF_{p,\pi_{G_0}}(G0) =\overline{\{\pi(f)\ :\ f\in L_1(G),\;\text{supp}(f)\subseteq G_0 \}{}}^{\|\cdot\|_{\mathcal{B}(E)}}\subseteq UPF_{p,\pi}(G).
\end{align*}
 \item\label{PROPEXREC2}
The linear restriction mapping $R_\pi : B_{p,\pi}\rightarrow B_{p,\pi_{G_0}}$ which is defined for $u\in B_{p,\pi}$, as $R_\pi(u) =u|_{G_0}$ is the dual map of $S_{\pi_{G_0}}$, and is a quotient
map.
\item\label{PROPEXREC2.5}
The extension map $E_\pi: B_{p,\pi_{G_0}}\rightarrow B_{p,\pi}$, defined via $E_\pi(u)=u^\circ$ is an isometric map.
\item\label{PROPEXREC3}
The restriction mapping $R : B_p(G)\rightarrow B_p(G_0)$ is a contraction.
\item\label{PROPEXREC4}
When $(\pi,E)$ is also a $p$-universal representation, we have the following contractive inclusions
\begin{align*}
PF_p(G_0)^*\subseteq B_{p,\pi_{G_0}}\subseteq B_p(G_0)\subseteq\mathcal{M}(A_p(G_0)).
\end{align*}
Under the assumption that $G_0$ is amenable, we have isometric identification below
\begin{align*}
PF_p(G_0)^* = B_{p,\pi_{G_0}}= B_p(G_0)=\mathcal{M}(A_p(G_0)).
\end{align*}
\end{enumerate}
\begin{proof}
\begin{enumerate}
\item
Through Lemma \ref{LEMMARESTRICTIONMAP}, the map $S_{\pi_{G_0}}$ is an isometric homomorphism with the range containing the dense space $\{\pi(f)\ :\ f\in L_1(G),\;\text{supp}(f)\subseteq G_0 \}$. So, the algebra $PF_{p,\pi_{G_0}}(G_0)$ and the subalgebra $\overline{\{\pi(f)\ :\ f\in L_1(G),\;\text{supp}(f)\subseteq G_0 \}{}}^{\|\cdot\|_{\mathcal{B}(E)}}$
of $UPF_p(G)$ are identified.
\item
Evidently, we have $R_\pi$ = $S^*_{\pi_{G_0}}$, therefore, $S_{\pi_{G_0}}$ is a quotient map.
\item
Before showing $E_\pi$ is an isometric map, it is needed to take notice of the fact that since $R_\pi$ is onto, then we have
\begin{align*}
B_{p,\pi_{G_0}}=PF_{p,\pi_{G_0}}(G_0)^*=PF_{p,\pi}(G)^*/PF_{p,\pi_{G_0}}^\bot=B_{p,\pi}/PF_{p,\pi_{G_0}}^\bot .
\end{align*}
Obviously, $E_\pi$ is contraction. Furthermore, we have
\begin{align*}
\|u\|=\|R_\pi(E_\pi(u))\|\leq \|E_\pi(u)\|\leq\|u\|.
\end{align*}
Moreover, one may be inclined to gain this conclusion through the following argument. Define the map $\mathcal{E}_\pi:PF_{p,\pi}(G)\rightarrow PF_{p,\pi_{G_0}}(G_0)$ via $\mathcal{E}_\pi(\pi(f))=\pi_{G_0}(f|_{G_0})$, for $f\in L_1(G)$. Then $\mathcal{E}_\pi$ is contraction and $(\mathcal{E}_\pi)^*=E_\pi$. So, $E_\pi$ is well-defined and contraction, as well.

\item
From the part \eqref{PROPEXREC2}, the restriction map from $B_p(G)$ onto $B_{p,\pi_{G_0}}$ is a contraction, and due to the Remark \ref{REMARKOFDUALOFPFP}-\eqref{REMARKOFDUALOFPFP4}, for every $(\rho, F)\in\text{Rep}_p(G_0)$ the identity map from $B_{p,\rho}$ into $B_p(G_0)$ is a contraction then we have the result of this part.
\item
Let $u\in PF_{p}(G_0)^*$. Then by the part \eqref{PROPEXREC2.5} we have $u^\circ\in PF_p(G)^*$, and $PF_p(G)^*\subseteq B_p(G)$, contractively, via Remark \ref{REMARKOFDUALOFPFP}-\eqref{REMARKOFDUALOFPFP6}. Since $R(u^\circ)=u$, it follows that $u\in B_{p,\pi_{G_0}}$, where $(\pi, E)$ is a $p$-universal representation of $G$.\\
For the case that $G_0$ is amenable, since through aforementioned remark, we have $PF_p(G_0)^*$ = $B_p(G_0)$ , so
$PF_p(G_0)^*=B_{p,\pi_{G_0}}=B_p(G_0)$.
\end{enumerate}
\end{proof}
\end{proposition}

Next proposition is the consequences of previous one, and is one of the applicable result in dealing with problems about $p$-analog of the Fourier-Stieltjes algebras.

\begin{proposition}\label{PROPEXTENSION}
Let $G$ be a locally compact group and $G_0$ be its open
subgroup. Then
\begin{enumerate}
\item\label{PROPEXTENSION1}
the extension mapping $E_{MM} : \mathcal{M}(A_p(G_0))\rightarrow \mathcal{M}(A_p(G))$, defined for
$u \in \mathcal{M}(A_p(G_0))$ via $E_{MM} (u) = u^\circ$
is an isometric map.
\item\label{PROPEXTENSION2}
for every $u\in B_p(G_0)$, we have $u^\circ\in \mathcal{M}(A_p(G))$, and the map $E_{BM} :
B_p(G_0)\rightarrow \mathcal{M}(A_p(G))$, with $u\mapsto u^\circ$, is a contraction.
\item\label{PROPEXTENSION3}
if $G_0$ is also an amenable subgroup, then for every $u \in B_p(G_0)$, we have $u^\circ \in B_p(G)$, and the associated extending map $E_{BB} : B_p(G_0)\rightarrow B_p(G)$ is an isometric one.
\end{enumerate}
\begin{proof}
\begin{enumerate}
\item
By the following relation for $u\in\mathcal{M}(A_p(G_0))$ and $v\in A_p(G)$
\begin{align*}
u^\circ\cdot v=(u\cdot v|_{G_0})^\circ,
\end{align*}
it can be concluded that $u^\circ\in\mathcal{M}(A_p(G))$, and obviously we have $\|u^\circ\|_{\mathcal{M}(A_p(G))}=\|u\|_{\mathcal{M}(A_p(G_0))}$
\item
This part can be concluded by the inclusions in Proposition \ref{PROPEXREC}-\eqref{PROPEXREC4} and the part \eqref{PROPEXTENSION1}.
\item
Since $G_0$ is amenable, then by Proposition \ref{PROPEXREC}-\eqref{PROPEXREC4} (or directly from Remark \ref{REMARKOFDUALOFPFP}-\eqref{REMARKOFDUALOFPFP6}) we have the result.
\end{enumerate}
\end{proof}
\end{proposition}

One of the interesting problems on the Fourier-Stieltjes type algebras is to study weighted homomorphism associated to a piecewise affine map as it has been considered in \cite{ILIESPRONK2005} and \cite{ILIE2013}. At this aim, it is crucial to be sure that such a homomorphism is well-defined. Precisely, answering to the question that the homomorphism $\Phi_\alpha :B_p(G)\rightarrow B_p(H)$, defined via $\Phi_\alpha(u)=(u\circ\alpha)^\circ$, for $u\in B_p(G)$ is well-defined or not, would be precious. Here, $\alpha :Y\subseteq H\rightarrow G$ is a piecewise affine map. So, we give some preliminaries here.\\
For a locally compact topological group $H$, let $\Omega_0(H)$ denote the ring of subsets which generated by open cosets of $H$. By \cite{ILIE2013} we have
\begin{align}\label{OPENCOSETRINGS}
\Omega_0(H)=\left\{Y\backslash\cup_{i=1}^nY_i \: :
\begin{array}{ll}
&Y\:\text{is an open coset of }\: H,\\
& Y_1,\ldots ,Y_n\:\text{open subcosets of infinite index in}\: Y\\
\end{array}\right\}.
\end{align}
Moreover, for a set $Y\subseteq H$, by $\text{Aff}(Y)$ we mean the smallest coset containing $Y$, and if $Y=Y_0\backslash\cup_{i=1}^nY_i\in\Omega_0(H)$, then $\text{Aff}(Y)=Y_0$.
Similarly, let us denote by  $\Omega_{\text{am-}0}(H)$ the ring of open cosets of open amenable subgroups of $H$. Now, we give the definition of a piecewise affine map.

\begin{definition}\label{DEFPIECEWISE}
Let $\alpha : Y\subseteq H\rightarrow G$ be a map.
\begin{enumerate}
\item
The map $\alpha$ is called an affine map on an open coset $Y$ of an open subgroup $H_0$, if
\begin{equation*}
\alpha(xy^{-1}z)=\alpha(x)\alpha(y)^{-1}\alpha(z),\qquad x,y,z\in Y,
\end{equation*}

\item
The map $\alpha$ is called a piecewise affine map if
\begin{enumerate}
\item
there are pairwise disjoint $ Y_i\in\Omega_0(H)$, for $ i=1,\ldots , n $, such that $Y=\cup_{i=1}^nY_i$,

\item there are affine maps $\displaystyle{\alpha_i : \text{Aff}(Y_i)\subseteq H\rightarrow G}$, for $ i=1,\ldots , n $, such that
\begin{equation*}
\alpha |_{Y_i}=\alpha_i |_{Y_i}.
\end{equation*}

\end{enumerate}
\end{enumerate}

\end{definition}

\begin{remark}\label{AFFFIINEREM}{\cite[Remark 2.2]{ILIE2004}}
If $Y=h_0H_0$ is an open coset of an open subgroup $H_0\subset H$, and $\alpha : Y\subseteq H\rightarrow G$ is an affine map, then there exists a group homomorphism $ \beta $ associated to $\alpha$ such that
\begin{align}\label{affine-homomorphism}
&\beta : H_0\subseteq H\rightarrow G,\quad\beta (h)= \alpha(h_0)^{-1}\alpha(h_0h),\quad h\in H_0.
\end{align}
\end{remark}
Next lemma is straightforward and we leave it without proof, and it will be utilized in Theorem \ref{THEOREMCONCLUSION}.

\begin{lemma}\label{LEMMAMAPMAP}
Let $G$ and $H$ are locally compact groups and $(\pi,E)\in\text{Rep}_p(G)$.
\begin{enumerate}
\item\label{LEMMAMAPMAP1}
For an element $x\in G$ let $L_x:B_p(G)\rightarrow B_p(G)$ be  the left translation mapping defined through $L_x(u)(y)=u(xy)$ for $y\in G$ and $u\in B_p(G)$. Then $L_x$ is an invertible isometric map.
\item\label{LEMMAMAPMAP2}
For a continuous homomorphism $\beta : H\rightarrow G$ the pair $(\pi\circ\beta, E)$ belongs to $\text{Rep}_p(H)$ and the homomorphism $\Phi_\beta :B_p(G)\rightarrow B_p(H)$ is well-defined contractive homomorphism.
\end{enumerate}
\end{lemma}

The following theorem is one our important results in this paper. Here, for a continuous piecewise affine map $\alpha :Y\subseteq H\rightarrow G$, we prove that the homomorphism $\Phi_\alpha:B_p(G)\rightarrow B_p(H)$, defined via
\begin{align*}
\Phi_\alpha(u)=\left\{
\begin{array}{lll}
u\circ\alpha &&\text{on}\; Y\\
0 && \text{o.w.}\; Y
\end{array}\right.
\end{align*}
is well-defined and we determine its bound.
\begin{theorem}\label{THEOREMCONCLUSION}
Let $G$ and $H$ be locally compact groups, and $\alpha :Y=\cup_{k=1}^nY_k\subseteq H\rightarrow G$ be a continuous piecewise affine map with disjoint $Y_k\in\Omega_{\text{am-}0}(H)$, for $k=1,\ldots,n$. Then  $u\in B_p(G)$ implies that $(u\circ \alpha)^\circ\in B_p(H)$, and consequently, the weighted homomorphism $\Phi_\alpha:B_p(G)\rightarrow B_p(H)$ is well-defined bounded homomorphism.
\begin{proof}
We divide our proof into two steps. 
\item[\textit{Step 1}:]
First, we let $\alpha: Y=y_0H_0\subseteq H\rightarrow G $ be a continuous affine map, and $\beta :H_0\rightarrow G$ be the homomorphism associated with $\alpha$, as it is explained in Remark \ref{AFFFIINEREM}, for an open amenable subgroup $H_0$ of $H$. As we initially explained in Lemma \ref{LEMMAMAPMAP}-\eqref{LEMMAMAPMAP2}, the map $u\mapsto u\circ\beta$ is an algebra homomorphism from $\Phi_\beta :B_p(G)\rightarrow B_p(H_0)$. For the element $y_0$ consider the  translation map $ L_{\alpha(y_0)}:  B_p(G)\rightarrow B_p(G)$, then by the following relation, and applying Proposition \ref{PROPEXTENSION}-\eqref{PROPEXTENSION3}, we have the result
\begin{align*}
(u\circ\alpha)^\circ=E_{BB}\circ \Phi_\beta\circ L_{\alpha(y_0)},\quad u\in B_p(G).
\end{align*}
where $E_{BB}:B_p(H_0)\rightarrow B_p(H)$, is the extension mapping. By the last relation it is obtained that the extension of the function $u\circ\alpha$ belongs to $B_p(H)$, and evidently $\Phi_\alpha$ is contractive as it is the combination of isometric and contractive maps.
\item[\textit{Step 2}:]
Now, let $\alpha: Y\subseteq H\rightarrow G $ be a continuous piecewise affine map, so by our assumption of amenability, and similar to the Definition \ref{DEFPIECEWISE}, there exist pairwise disjoint sets $Y_k\in\Omega_{\text{am-}0}(H)$, for $k=1,\ldots,n$ with $n\in\mathbb{N}$, and affine maps $\alpha_k:\text{Aff}(Y_k)\subseteq H\rightarrow G$ such that $Y=\cup_{k=1}^n Y_k$, and $\alpha_k|_{Y_k}=\alpha|_{Y_k}$. By previous step, we know that $(u\circ\alpha_k)^\circ\in B_p(H)$, and since 
\begin{align*}
(u\circ\alpha_k)^\circ=\sum_{k=1}^n(u\circ\alpha_k)^\circ\cdot\chi_{Y_k},
\end{align*}
we have the result via Corollary \ref{COROIDEMCOSETRING}, and the fact that $B_p(H)$ is a Banach algebra. Moreover, we have
\begin{align*}
\|(u\circ\alpha)^\circ\|\leq\|u\|\sum_{k=1}^n 2^{m_{Y_k}},
\end{align*}
where the number $m_{Y_k}$ is as it is described in Corollary \ref{COROIDEMCOSETRING}. So, we have $\|\Phi_\alpha\|\leq \sum_{k=1}^n 2^{m_{Y_k}}$. 
\end{proof}
\end{theorem}

\begin{remark}\label{REMARKEXTENSION}
For an open amenable subgroup $G_0$ of the locally compact group $G$, by Proposition \ref{PROPEXTENSION}-\eqref{PROPEXTENSION3}, we can say that the space $B_p(G_0)$ is the space of functions which are restriction of functions in $B_p(G)$, those are equal to zero outside of $G_0$. Therefore, while we are working on the $p$-analog of the Fourier-Stieltjes algebras, we may assume that the $p$-universal representation of an open amenable subgroup $G_0$ of $G$ is the restriction of the $p$-universal representation of $G$ to $G_0$.
\end{remark}

\end{document}